\documentclass[a4paper]{amsart}
\usepackage[british]{babel}
\usepackage[utf8]{inputenc}
\usepackage{csquotes}

\usepackage{amsmath,amssymb,amsthm}
\usepackage{tikz}
\usepackage{tikz-cd}
\usepackage[style=alphabetic,sorting=anyt,sortcites=true,backend=biber,maxbibnames=5,maxcitenames=5, giveninits = true]{biblatex}
\usepackage{float}
\usepackage{hyperref}

\usepackage{algorithm}
\usepackage{algpseudocode}
\usepackage{cleveref}

\usepackage{enumerate}

\makeatletter
\newcounter{algorithmicH}
\let\oldalgorithmic\algorithmic
\renewcommand{\algorithmic}{%
	\stepcounter{algorithmicH}
	\oldalgorithmic}
\renewcommand{\theHALG@line}{ALG@line.\thealgorithmicH.\arabic{ALG@line}}
\makeatother

\newcommand{\AlgFail}{\texttt{false}}

\theoremstyle{plain}
\newtheorem{theorem}{Theorem}[section]
\newtheorem{lemma}[theorem]{Lemma}
\newtheorem{proposition}[theorem]{Proposition}
\newtheorem{corollary}[theorem]{Corollary}

\theoremstyle{definition}

\theoremstyle{remark}

\DeclareMathOperator{\End}{End}
\DeclareMathOperator{\Aut}{Aut}

\DeclareMathOperator{\coker}{coker}
\DeclareMathOperator{\im}{im}

\DeclareMathOperator{\Fitt}{Fitt}
\DeclareMathOperator{\id}{id}

\newcommand{\R}{\mathfrak{R}}

\DeclareMathOperator{\Coin}{Coin}

\newcommand{\bvarphi}{\bar{\varphi}}
\newcommand{\bpsi}{\bar{\psi}}

\title[{Algorithms for twisted conj. classes of polycyclic-by-finite groups}]{Algorithms for twisted conjugacy classes of polycyclic-by-finite groups}
\author[K. Dekimpe]{Karel Dekimpe}
\author[S. Tertooy]{Sam Tertooy}
\address{KU Leuven Campus Kulak Kortrijk\\
	E. Sabbelaan 53\\
	8500 Kortrijk\\
	Belgium}
\email{Karel.Dekimpe@kuleuven.be}
\email{Sam.Tertooy@kuleuven.be}
\thanks{Research supported  by long term structural funding -- Methusalem grant of the Flemish Government.}


\begin{document}

		\begin{abstract}
		We construct two practical algorithms for twisted conjugacy classes of polycyclic groups. The first algorithm determines whether two elements of a group are twisted conjugate for two given endomorphisms, under the condition that their Reidemeister coincidence number is finite. The second algorithm determines representatives of the Reidemeister coincidence classes of two endomorphisms if their Reidemeister coincidence number is finite, or returns ``false'' if this number is infinite. We also discuss a theoretical extension of these algorithms to polycyclic-by-finite groups.
	\end{abstract}

	\maketitle

	\begin{center}
	This is an Accepted Manuscript of an article published by Elsevier in Topology and its Applications on 15 Apr 2021, available online:  \href{https://doi.org/10.1016/j.topol.2020.107565}{https://doi.org/10.1016/j.topol.2020.107565}.
\end{center}
\section{Introduction}
\label{sec:intro}
	Let \(G\) and \(H\) be groups and let \(\varphi\), \(\psi: H \to G\) be group homomorphisms. The coincidence group \(\Coin(\varphi,\psi)\) of the pair \((\varphi,\psi)\) is the subgroup of \(H\) defined by
	\begin{equation*}
	\Coin(\varphi,\psi) := \{ h \in H \mid \varphi(h) = \psi(h) \}.
	\end{equation*}
	Define an equivalence relation \(\sim_{\varphi,\psi}\) on \(G\) by
	\begin{equation*}
	\forall g_1,g_2 \in G: g_1 \sim_{\varphi,\psi} g_2 \iff \exists h \in H: g_1 = \psi(h)g_2\varphi(h)^{-1}.
	\end{equation*}
	The equivalence classes \([g]_{\varphi,\psi}\) are called the \emph{Reidemeister (coincidence) classes} of the pair \((\varphi,\psi)\) or the \emph{\((\varphi,\psi)\)-twisted conjugacy classes}. The set of Reidemeister classes is denoted by \(\R(\varphi,\psi)\). The \emph{Reidemeister (coincidence) number} \(R(\varphi,\psi)\) is the cardinality of \(\R(\varphi,\psi)\) and is therefore always a positive integer or infinity.

This equivalence relation originates in topological coincidence theory, see \cite{gonc05-1} for a survey. One of the aims of coincidence theory is, given two continuous maps \(f\), \(g: X \to Y\) between topological spaces \(X\), \(Y\), to calculate the number
\begin{equation*}
MC(f,g) := \min_{f' \simeq f, g' \simeq g} \# \{ x \in X \mid f'(x) = g'(x) \},
\end{equation*}
i.e. the least number of coincidence points among any pair of maps \((f',g')\), with \(f'\) in the homotopy class of \(f\) and \(g'\) in the homotopy class of \(g\). The \emph{Nielsen coincidence number} \(N(f,g)\), defined as the number of essential coincidence classes of the pair \((f,g)\), is a lower bound for \(MC(f,g)\). The \emph{Reidemeister coincidence number} \(R(f,g)\), defined as the number of coincidence classes (essential or otherwise) of the pair \((f,g)\), is an upper bound for the Nielsen coincidence number. While the Nielsen number will always be finite and can be zero, the Reidemeister number is either positive or infinite. In general, Nielsen numbers are quite difficult to compute, whereas Reidemeister numbers are much easier to calculate. The Reidemeister coincidence  number \(R(f,g)\) of continuous maps \(f\) and \(g\) equals the Reidemeister coincidence number \(R(f_*,g_*)\) of the induced group homomorphisms \(f_*\), \(g_*:\pi_1(X) \to \pi_1(Y)\) between the fundamental groups of \(X\) and \(Y\). 

If \(f\), \(g: M \to M\) are continuous self-maps of an orientable infra-nilmanifold or an infra-solvmanifold \(M\) of type (R), or if \(g = \id_M\) and \(f\) is a continuous self-map of any infra-solvmanifold \(M\), then the Nielsen coincidence number \(N(f,g)\) equals the Reidemeister coincidence number \(R(f,g)\) if the latter is finite, see \cite{dp11-1,fl15-1,dv20-1}. The fundamental group of an infra-solvmanifold is a (torsion-free) polycyclic-by-finite group, and conversely, every torsion-free polycyclic-by-finite group is the fundamental group of some infra-solvmanifold \cite{baue04-1}.

In \cite{su08-1}, an authentication scheme is proposed that relies on the ``apparent hardness of the twisted conjugacy problem'', i.e. given \(g_1 \sim_{\varphi,\psi} g_2\), it is assumed to be difficult to calculate some \(h\) such that \(g_1 = \psi(h)g_2\varphi(h)^{-1}\). Polycyclic groups have been suggested as the platform groups for various cryptosystems, including this authentication scheme \cite{gk16-1}.

The main goal of this paper is to construct two algorithms for endomorphisms of polycyclic-by-finite groups, which will be practical when applied to polycyclic groups. The first algorithm, which we will call \textsc{RepTwistConj} (short for \emph{Representative for Twisted Conjugation}), takes as input two endomorphisms \(\varphi\), \(\psi: G \to G\) with finite Reidemeister number \(R(\varphi,\psi)\) and two elements \(g_1\), \(g_2\) of a polycyclic-by-finite group \(G\), and returns the following output:
\begin{itemize}
	\item if \(g_1 \sim_{\varphi,\psi} g_2\): an element \(h \in G\) such that \(g_1 = \psi(h)g_2\varphi(h)^{-1}\),
	\item if \(g_1 \not\sim_{\varphi,\psi} g_2\): ``{\AlgFail}''. 
\end{itemize}
The second algorithm, which we will call \textsc{RepsReidClasses} (short for \emph{Representatives of Reidemeister Classes}), takes as input two endomorphisms \(\varphi\), \(\psi: G \to G\) of a polycyclic-by-finite group \(G\) and returns the following output:
\begin{itemize}
	\item if \(R(\varphi,\psi) < \infty\): a finite subset \(\{g_1, \ldots, g_n\} \subseteq G\) for which \(g_i \not\sim_{\varphi,\psi} g_j\) when \(i \neq j\) and \(\R(\varphi,\psi) = \{[g_1]_{\varphi,\psi}, \ldots, [g_n]_{\varphi,\psi}\}\),
	\item if \(R(\varphi,\psi) = \infty\): ``{\AlgFail}''.
\end{itemize}
Together, these algorithms will determine whether or not the Reidemeister coincidence number is finite, and if it is, they will completely determine the Reidemeister coincidence classes. In particular, this allows us to calculate Reidemeister coincidence numbers of polycyclic-by-finite groups, and thus Reidemeister coincidence numbers of infra-solvmanifolds as well. Moreover, these algorithms demonstrate that if a polycyclic group is used as platform group for the authentication scheme from \cite{su08-1}, then the endomorphisms should be picked such that they have infinite Reidemeister coincidence number.

In this paper, we will assume that the group \(G\) belongs to a class of groups suitable for computation with homomorphisms, in the sense that we can easily do the following:
\begin{itemize}
	\item calculate the kernel and image of homomorphisms between groups in this class,
	\item calculate the image or a preimage of a group element under the above homomorphisms.
\end{itemize}
The class of polycyclic groups satisfies these criteria, and practical algorithms to do the above can found in the GAP-package \texttt{polycyclic} \cite{gap18-5} and the references contained in its package manual.

\section{Preliminaries}
Throughout this paper, we will use the notation \(\iota_x\) to describe the inner automorphism \(G \to G: g \mapsto xgx^{-1}\).
\begin{lemma}
	Let \(G\) be a group, \(\varphi,\psi \in \End(G)\) and \(g_1,g_2 \in G\). For any \(x \in G\), we have that
	\begin{equation*}
	g_1 \sim_{\varphi,\psi} g_2 \iff g_1x^{-1} \sim_{\iota_x\varphi,\psi} g_2x^{-1},
	\end{equation*}
	and moreover
	\begin{equation*}
	\{h \in G \mid g_1 = \psi(h)g_2\varphi(h)^{-1}\} = \{h \in G \mid g_1x^{-1} = \psi(h)g_2x^{-1}(\iota_x\varphi)(h)^{-1}\}.
	\end{equation*}
\end{lemma}
\begin{proof}
	For any \(h \in G\), we have that
	\begin{align*}
	g_1 = \psi(h)g_2\varphi(h)^{-1} &\iff g_1x^{-1} = \psi(h)g_2x^{-1}x\varphi(h)^{-1}x^{-1}\\
	&\iff g_1x^{-1} = \psi(h)g_2x^{-1}(\iota_x\varphi)(h)^{-1}.\qedhere
	\end{align*}
	
\end{proof}
By taking \(x = g_2\) in the above lemma, we obtain the following corollary.
\begin{corollary}
	\label{cor:twistconjtoidentity}
	 Let \(g_1,g_2 \in G\) and \(\varphi,\psi \in \End(G)\). Then \(g_1 \sim_{\varphi,\psi} g_2\) if and only if \(g_1g_2^{-1} \sim_{\iota_{g_2}\varphi,\psi} 1\).
\end{corollary}
Thus, it suffices to solve the twisted conjugacy problem in the case where one of the elements is the identity. This does, however, involve composing one of the endomorphisms with an inner automorphism. The following corollary shows that this does not impact the finiteness of the Reidemeister coincidence number of the endomorphisms.
\begin{corollary}
	\label{cor:psiisbijection}
	Let \(g \in G\) and let \(\varphi,\psi \in \End(G)\). Then the map \(\mu_g:  \R(\iota_g\varphi,\psi) \to \R(\varphi,\psi): [x]_{\iota_g\varphi,\psi} \mapsto [xg]_{\varphi,\psi}\) is a bijection, and therefore \(R(\varphi,\psi) = R(\iota_g\varphi,\psi)\).
\end{corollary}
Therefore, should we have an algorithm \textsc{RepTwistConjToId}(\(\varphi,\psi,g)\) that takes as input two endomorphisms \(\varphi\), \(\psi\) with finite Reidemeister number \(R(\varphi,\psi)\) and an element \(g\), and returns the following output:
\begin{itemize}
	\item if \(g \sim_{\varphi,\psi} 1\): an element \(h \in G\) such that \(g = \psi(h)\varphi(h)^{-1}\),
	\item if \(g \not\sim_{\varphi,\psi} 1\): ``{\AlgFail}'',
\end{itemize}
then we may construct the algorithm \textsc{RepTwistConj} as in \Cref{alg:RepTwistConj}. 
\begin{algorithm}
	\caption{Determining \(h\) such that \(g_1 = \psi(h)g_2 \varphi(h)^{-1}\)}
	\label{alg:RepTwistConj}
	\begin{algorithmic}[1]
		\Function{RepTwistConj}{$\varphi,\psi,g_1,g_2$}
		\State \Return \( \textsc{RepTwistConjToId}(\iota_{g_2}\varphi,\psi,g_1g_2^{-1})\)
		\EndFunction
	\end{algorithmic}
\end{algorithm}

The following theorem will be crucial in constructing both \textsc{RepTwistConjToId} and \textsc{RepsReidClasses} for polycyclic and polycyclic-by-finite groups.
\begin{theorem}[see {\cite[\textsection2]{kl07-1}}]
	\label{thm:diagramExactSequence}
	Let \(G\) be group, let \(N\) be a normal subgroup of \(G\) and let \(\varphi,\psi \in \End(G)\) such that \(\varphi(N) \subseteq N\) and \(\psi(N) \subseteq N\). We denote the restrictions of \(\varphi\) and \(\psi\) to \(N\) by \(\varphi|_N\) and \(\psi|_N\), and the induced endomorphisms on \(G/N\) by \(\bvarphi\) and \(\bpsi\). We then get the following commutative diagram with exact rows:
	\begin{center}
		\begin{tikzcd}[row sep=large,column sep=large]
		1 \arrow{r} & N \arrow{r}{i}\arrow[shift left]{d}[shift left=13]{\varphi|_N}\arrow[shift right]{d}[swap]{\psi|_N} & G                 \arrow{r}{p}\arrow[shift left]{d}{\varphi}\arrow[shift right]{d}[swap]{\psi}      & G/N           \arrow{r}\arrow[shift left]{d}{\bvarphi}\arrow[shift right]{d}[swap]{\bpsi} & 1 \\
		1 \arrow{r} & N \arrow{r}{i}                & G             \arrow{r}{p}      & G/N        \arrow{r} & 1
		\end{tikzcd}
		
	\end{center}
	This diagram induces the following exact sequence of pointed sets:

	\begin{center}
	\begin{tikzcd}[row sep=large]
	1 \arrow{r} & \Coin(\varphi|_N,\psi|_N) \arrow{r}{i} & \Coin(\varphi,\psi)  \arrow{r}{p}      & \Coin(\bvarphi,\bpsi)  \arrow[out=0, in=180, looseness=2,overlay]{dll}[swap]{\delta} &  \\
	& \R(\varphi|_N,\psi|_N) \arrow{r}{\hat{\imath}}    & \R(\varphi,\psi)  \arrow{r}{\hat{p}}      & \R(\bvarphi,\bpsi)        \arrow{r} & 1
	\end{tikzcd}
\end{center}
	where all maps are evident except \(\delta\), which is defined as \(\delta(\bar{g}) = [\psi(g)\varphi(g)^{-1}]_{\varphi|_N,\psi|_N}\). 
\end{theorem}
The corollary below is a straightforward generalisation of statements (1) and (2) in \cite[Lemma 1.1]{gw09-2}. 
\begin{corollary}
	\label{cor:2props}
		Consider the situation from \Cref{thm:diagramExactSequence}. We obtain the following properties:
		\begin{enumerate}[(1)]
		\item \(R(\varphi,\psi) \geq R(\bvarphi,\bpsi)\),
		\item if \(\# \Coin(\bvarphi,\bpsi) < \infty\) and \(R(\varphi,\psi) < \infty\), then \(R(\varphi|_N,\psi|_N) < \infty\).
	\end{enumerate}
\end{corollary}

\section{Reduction to normal subgroup and quotient}
\label{sec:reduction}

It is possible to reduce the twisted conjugacy problem on a group to the twisted conjugacy problem on a well-chosen normal subgroup and on the quotient by that subgroup.

\begin{theorem}
	\label{thm:reduceEquivToIdentityToNormalSubgroup}
	Consider the situation from \Cref{thm:diagramExactSequence}. Let \(g \in G\). If \(\bar{g} \sim_{\bvarphi,\bpsi} \bar{1}\), then there exists an \(n \in N\) such that \(n \sim_{\varphi,\psi} g\) and
	\begin{equation*}
	g \sim_{\varphi,\psi} 1 \iff \exists \bar{h} \in \Coin(\bvarphi,\bpsi): \psi(h)^{-1}n\varphi(h) \sim_{\varphi|_N,\psi|_N} 1,
	\end{equation*}
	where \(\bar{g} := p(g)\) and \(h\) is any element of \(p^{-1}(\bar{h})\).
\end{theorem}

\begin{proof}
	If \(\bar{g} \sim_{\bvarphi,\bpsi} \bar{1}\), then there exists some \(\bar{k} \in G/N\) such that
	\begin{equation*}
	\bar{g} = \bpsi(\bar{k})\bvarphi(\bar{k})^{-1} \iff \bpsi(\bar{k})^{-1}\bar{g}\bvarphi(\bar{k}) =   \bar{1}.
	\end{equation*}
	Let \(k \in G\) be any preimage of \(\bar{k}\), then \(n := \psi(k)^{-1}g\varphi(k)\) is an element of \(N\) and clearly \(n \sim_{\varphi,\psi} g\). Now, using the exact sequence from \Cref{thm:diagramExactSequence}, we find that
	\begin{align*}
	[g]_{\varphi,\psi} = [1]_{\varphi,\psi} &\iff [n]_{\varphi,\psi} = [1]_{\varphi,\psi}\\
	&\iff \hat{\imath}([n]_{\varphi|_N,\psi|_N}) = [1]_{\varphi,\psi}\\
	&\iff \exists \bar{h} \in \Coin(\bvarphi,\bpsi): [n]_{\varphi|_N,\psi|_N} = [\psi(h)\varphi(h)^{-1}]_{\varphi|_N,\psi|_N}\\
	&\iff \exists \bar{h} \in \Coin(\bvarphi,\bpsi): [\psi(h)^{-1}n\varphi(h)]_{\varphi|_N,\psi|_N} = [1]_{\varphi|_N,\psi|_N},
	\end{align*}
	where we used the normality of \(N\) to obtain the last equivalence.
\end{proof}
Thus, we can construct \Cref{alg:RepTwiConIdNormal}, called \textsc{RepTwistConjToIdByNormal}, which reduces the twisted conjugacy problem on \(G\) to the twisted conjugacy problem on a normal subgroup \(N\) and on the quotient \(G/N\). In order for this algorithm to work, we require \(4\) conditions on the endomorphisms \(\varphi, \psi\) and the normal subgroup \(N\) given as input:
\begin{enumerate}[(i)]\label{enum:conds}
	\item \(\varphi(N) \subseteq N\) and \(\psi(N) \subseteq N\), such that \(\bvarphi\), \(\bpsi\), \(\varphi|_N\) and \(\psi|_N\) are well-defined,
	\item \(\Coin(\bvarphi,\bpsi)\) must be finite and be easily computable, because \cref{alg2l:forloop} iterates over all elements of this group.
	\item \textsc{RepTwistConjToId} is implemented for input \(\bvarphi,\bpsi\), because \cref{alg2l:calltoquotient} calls this,
	\item \textsc{RepTwistConjToId} is implemented for input \(\varphi|_N,\psi|_N\), because \cref{alg2l:calltosubgroup} calls this.
\end{enumerate}
Note that we currently do not require that \(R(\varphi,\psi) < \infty\). 
\begin{algorithm}
	\caption{Determining \(h\) such that \(g = \psi(h) \varphi(h)^{-1}\)}
	\label{alg:RepTwiConIdNormal}
	\begin{algorithmic}[1]
		\Function{RepTwistConjToIdByNormal}{$\varphi,\psi,g,N$}
		\State \(p :=\) projection \(G \to G/N\)
		\State \(\bar{k} := \textsc{RepTwistConjToId}(\bvarphi,\bpsi,p(g))\) \label{alg2l:calltoquotient}
		\If{\(\bar{k} = \AlgFail\)}
		\State \Return \AlgFail
		\EndIf
		\State \(k :=\) any element in \(p^{-1}(\bar{k})\)
		\State \(n := \psi(k)^{-1}g\varphi(k)\)
		\For{\(\bar{h} \in \Coin(\bvarphi,\bpsi)\)} \label{alg2l:forloop}
		\State \(h :=\) any element in \(p^{-1}(\bar{h})\)
		\State \(l := \textsc{RepTwistConjToId}(\varphi|_N,\psi|_N,\psi(h)^{-1}n\varphi(h))\) \label{alg2l:calltosubgroup}
		\If{\(l \neq \AlgFail\)}
		\State \Return \(khl\)
		\EndIf
		\EndFor
		\State \Return \AlgFail
		\EndFunction
	\end{algorithmic}
\end{algorithm}

Making use of the exact sequence from \Cref{thm:diagramExactSequence}, we may describe the set of Reidemeister classes \(\R(\varphi,\psi)\) in terms of Reidemeister classes of a well-chosen normal subgroup and of the quotient by that subgroup.
\begin{theorem}
	\label{thm:reidclassesunion}
	Consider the situation from \Cref{thm:diagramExactSequence}. The set of Reidemeister classes of the pair \((\varphi,\psi)\) is given by
	\begin{equation*}
	\R(\varphi,\psi) = \bigsqcup_{[\bar{g}]_{\bvarphi,\bpsi} \in \R(\bvarphi,\bpsi)} (\mu_g \circ \hat{\imath}_g)(\R(\iota_g\varphi|_N,\psi|_N)),
	\end{equation*}
	where \(\hat{\imath}_g\) is the map
	\begin{equation*}
	\hat{\imath}_g: \R(\iota_g\varphi|_N,\psi|_N) \to  \R(\iota_g\varphi,\psi): [x]_{\iota_g\varphi|_N,\psi|_N} \to [x]_{\iota_g\varphi,\psi}
	\end{equation*}
	and \(\mu_g\) is the map from \Cref{cor:psiisbijection}.
\end{theorem}
\begin{proof}
	From the surjectivity of \(\hat{p}\), we have that
		\begin{equation}
		\label{eq:rphidisjointunion}
	\R(\varphi,\psi) = \bigsqcup_{[\bar{g}]_{\bvarphi,\bpsi} \in \R(\bvarphi,\bpsi)}\hat{p}^{-1}([\bar{g}]_{\bvarphi,\bpsi}).
	\end{equation}
	Let \(\hat{p}_g\) be the map
	\begin{equation*}
	\hat{p}_g: \R(\iota_g\varphi,\psi) \to  \R(\iota_{\bar{g}}\bvarphi,\bpsi): [x]_{\iota_g\varphi,\psi} \to [\bar{x}]_{\iota_{\bar{g}}\bvarphi,\bpsi},
	\end{equation*}
	then by \Cref{cor:psiisbijection} and the exact sequence from \Cref{thm:diagramExactSequence} we obtain that \begin{equation}
	\label{eq:pinverse}
	\hat{p}^{-1}([\bar{g}]_{\bvarphi,\bpsi})   = \mu_g( \hat{p}_{g}^{-1}([\bar{1}]_{\iota_{\bar{g}}\bvarphi,\bpsi}) )= (\mu_g \circ \hat{\imath}_g)(\R(\iota_g\varphi|_N,\psi|_N)).
	\end{equation}
	The result now follows by combining \eqref{eq:rphidisjointunion} and \eqref{eq:pinverse}.
\end{proof}
Similar to the previous algorithm, we can construct \Cref{alg:RepReidClassNormal}. This time, we require \(5\) conditions on the endomorphisms \(\varphi, \psi\) and the normal subgroup \(N\) given as input in order for this algorithm to work as intended:
\begin{enumerate}[(i)]\label{enum:conds2}
	\item \(\varphi(N) \subseteq N\) and \(\psi(N) \subseteq N\), such that \(\bvarphi\), \(\bpsi\), \(\varphi|_N\) and \(\psi|_N\) are well-defined,
	\item \textsc{RepsReidClasses} is implemented for input \(\bvarphi,\bpsi\), because \cref{alg3l:calltoquotient} calls this,
	\item \textsc{RepsReidClasses} is implemented for input \(\iota_g\varphi|_N,\psi|_N\), because \cref{alg3l:calltosubgroup} calls this,
	\item If \(R(\bvarphi,\bpsi) < \infty\) and \(R(\iota_g\varphi|_N,\psi|_N) = \infty\) for some \(g \in G\), then \(R(\varphi,\psi) = \infty\), because \cref{alg3l:assumption} makes this assumption,
	\item \textsc{RepTwistConj} is implemented for input \(\iota_g\varphi,\psi\), because \cref{alg3l:calltotwistconj} calls this.
\end{enumerate}

\begin{algorithm}
	\caption{Determining representatives of \(\R(\varphi,\psi)\)}
	\label{alg:RepReidClassNormal}
	\begin{algorithmic}[1]
		\Function{RepsReidClassesByNormal}{$\varphi,\psi,N$}
		\State \(p :=\) projection \(G \to G/N\)
		\State \(\R(\bvarphi,\bpsi) := \textsc{RepsReidClasses}(\bvarphi,\bpsi)\)\label{alg3l:calltoquotient}
		\If{\(\R(\bvarphi,\bpsi)= \AlgFail\)}
		\State \Return \AlgFail
		\EndIf
		\State \(\R := \varnothing\)
		\For{\(\bar{g} \in \R(\bvarphi,\bpsi)\)}
		\State \(g :=\) any element in \(p^{-1}(\bar{g})\)
		\State \(\R(\iota_g\varphi|_N,\psi|_N) := \textsc{RepsReidClasses}(\iota_g\varphi|_N,\psi|_N)\) \label{alg3l:calltosubgroup}
		\If{\(\R(\iota_g\varphi|_N,\psi|_N)= \AlgFail\)}
		\State \Return  \AlgFail\label{alg3l:assumption}
		\EndIf
		\State \(\hat{\imath}_g\R := \varnothing\)
		\For{\(h \in \R(\iota_g\varphi|_N,\psi|_N)\)}
		\If{\(\forall k \in \hat{\imath}_g \R: \textsc{RepTwistConj}(\iota_g\varphi,\psi,h,k) = \AlgFail\)}  \label{alg3l:calltotwistconj}
		\State \(\hat{\imath}_g \R := \hat{\imath}_g \R \cup \{h\}\)
		\EndIf
		\EndFor
		\State \(\R := \R \cup \mu_g(\hat{\imath}_g \R)\)
		\EndFor
		\State \Return \(\R\)
		\EndFunction
	\end{algorithmic}
\end{algorithm}

\section{Abelian Groups}
If the group \(G\) is abelian, the set of Reidemeister classes can actually be interpreted as a quotient group of \(G\).
\begin{theorem}
	\label{thm:reidclassesabelian}
Let \(G\) be an abelian group and \(\varphi,\psi \in \End(G)\). Then \(\R(\varphi,\psi) = \coker(\psi-\varphi)\).
\end{theorem}
\begin{proof}
	Let \(g_1,g_2 \in G\). Then
	\begin{align*}
	g_1 \sim_{\varphi,\psi} g_2 &\iff \exists h \in G: g_1= \psi(h)+g_2 - \varphi(h) \\
	&\iff \exists h \in G:  g_1-g_2 = (\psi-\varphi)(h)\\
	&\iff g_1+\im(\psi-\varphi) = g_2 + \im(\psi-\varphi). \qedhere
	\end{align*}
\end{proof}
Thus, we can define \textsc{RepTwistConjToId} and \textsc{RepsReidClasses} for finitely generated, abelian groups as in \Cref{alg:RepTwiConIdAbelian,alg:RepReidClassAbelian}.

\begin{algorithm}
	\caption{Determining \(h\) such that \(g = \psi(h) \varphi(h)^{-1}\) if \(G\) is abelian}
	\label{alg:RepTwiConIdAbelian}
	\begin{algorithmic}[1]
		\Function{RepTwistConjToId}{$\varphi,\psi,g$}
		\If{\(g \in \im(\psi-\varphi)\)}
			\State \(h := \) any element in \((\psi-\varphi)^{-1}(g)\)
			\State \Return \(h\)
		\EndIf
		\State \Return \AlgFail
		\EndFunction
	\end{algorithmic}
\end{algorithm}

\begin{algorithm}
	\caption{Determining representatives of \(\R(\varphi,\psi)\) if \(G\) is abelian}
	\label{alg:RepReidClassAbelian}
	\begin{algorithmic}[1]
		\Function{RepsReidClasses}{$\varphi,\psi$}
		\If{\([G:\im(\psi-\varphi)] = \infty\)}
			\State \Return \AlgFail
		\EndIf
		\State \(\R:= \varnothing\)
		\State \(p :=\) projection \(G \to G/\im(\psi-\varphi)\)
		\For{\(\bar{g} \in G/\im(\psi-\varphi)\)}
		\State \(g :=\) any element in \(p^{-1}(\bar{g})\)
		\State \(\R := \R \cup \{g\}\)
		\EndFor
		\State \Return \(\R\)
		\EndFunction
	\end{algorithmic}
\end{algorithm}

The following proposition and corollary will be necessary when dealing with abelian quotients of polycyclic groups.
\begin{proposition}
	\label{prop:rankkercoker}
	Let \(G\) be a finitely generated, abelian group and let \(\varphi \in \End(G)\). Then the Hirsch length of the kernel of \(\varphi\) equals the Hirsch length of the cokernel of \(\varphi\).
\end{proposition}
\begin{proof}
	It is well known that for any polycyclic group with normal subgroup \(N\), \(h(G) = h(N) + h(G/N)\). Since \(\im(\varphi) \cong G/\ker(\varphi)\) and \(\coker(\varphi) = G/\im(\varphi)\), we obtain
	\begin{equation*}
	h(\ker(\varphi)) + h(\im(\varphi)) = h(G) = h(\im(\varphi)) + h( \coker(\varphi)).
	\end{equation*}
	Subtracting \(h(\im(\varphi))\) from both sides gives us the desired result.
\end{proof}

\begin{corollary}
	\label{cor:abelianRphifiniteiffFixphifinite}
	Let \(G\) be a finitely generated, abelian group and let \(\varphi,\psi \in \End(G)\). Then \(R(\varphi,\psi)\) is finite if and only if \(\Coin(\varphi,\psi)\) is finite.
\end{corollary}
\begin{proof}
	Note that \(\R(\varphi,\psi) = \coker(\psi-\varphi)\) (see \Cref{thm:reidclassesabelian}) and that \(\Coin(\varphi,\psi) = \ker(\psi-\varphi)\). By \Cref{prop:rankkercoker}, if either of these groups has Hirsch length \(0\), then so does the other.
\end{proof}

\section{Polycyclic groups}
One way to define a polycyclic group, is to state that all of its subgroups are finitely generated and that its derived series terminates at the trivial subgroup. This derived series will be exceptionally useful in the context of twisted conjugacy, as every group in this series is fully invariant and the factors are finitely generated, abelian groups.

\begin{proposition}
	\label{prop:recursiveReidNrsPolycyclic}
	Consider the situation from \Cref{thm:diagramExactSequence}, where \(G\) and \(N\) are chosen in such way that \(G/N\) is a finitely generated, abelian group. If \(R(\varphi,\psi)\) is finite, then so are \(\#\Coin(\bvarphi,\bpsi)\), \(R(\bvarphi,\bpsi)\) and \(R(\varphi|_N,\psi|_N)\).
\end{proposition}
\begin{proof}
	If \(R(\varphi,\psi) < \infty\), then by \Cref{cor:2props}(1) \(R(\bvarphi,\bpsi) < \infty\) and thus \Cref{cor:abelianRphifiniteiffFixphifinite} gives us that \(\# \Coin(\bvarphi,\bpsi) < \infty\). Finally, by \Cref{cor:2props}(2) \(R(\varphi|_N,\psi|_N)\) is finite as well.
\end{proof}
 \Cref{alg:RepTwiConIdPolycyclic} provides an implementation of \textsc{RepTwistConjToId} for polycyclic groups of derived length at least \(2\), under the restriction that the pair of endomorphisms given as input has finite Reidemeister number.
 
\begin{theorem}
	\label{thm:algworks}
	Let \(G\) be a polycyclic group of derived length at least \(2\) and let \(\varphi,\psi \in \End(G)\) such that  \(R(\varphi,\psi) < \infty\). Then \(\varphi\), \(\psi\) and \(G'\) satisfy the conditions necessary to apply \Cref{alg:RepTwiConIdNormal}.
\end{theorem}
\begin{proof}We prove this condition by condition. 
	\begin{enumerate}[(i)]
		\item This condition is satisfied because the derived subgroup \(G'\) is fully invariant.
		\item Since \(R(\varphi,\psi) < \infty\) and \(G/G'\) is finitely generated and abelian, \Cref{prop:recursiveReidNrsPolycyclic} gives us that \(\Coin(\bvarphi,\bpsi)\) is finite. As \(\Coin(\bvarphi,\bpsi) = \ker(\bpsi - \bvarphi) \subseteq G/G'\) it can be computed effectively.
		\item \Cref{alg:RepTwiConIdAbelian} provides an implementation for endomorphisms of \(G/G'\). 
		\item  We prove this by induction on the derived length \(n\) of \(G\). If \(n =2\), then \(G'\) is abelian, hence \Cref{alg:RepTwiConIdAbelian} provides an implementation for endomorphisms of \(G'\). Now assume that \(G\) has derived length \(n\) and that this theorem holds if the derived length is at most \(n-1\). By \Cref{prop:recursiveReidNrsPolycyclic} and the induction hypothesis, \(\varphi|_{G'}\), \(\psi|_{G'}\) and \(G''\) satisfy conditions (i) - (iv), thus \Cref{alg:RepTwiConIdPolycyclic} provides an implementation.\qedhere
	\end{enumerate}
\end{proof}

\begin{algorithm}
	\caption{Determining \(h\) such that \(g = \psi(h) \varphi(h)^{-1}\) if \(G\) is polycyclic}
	\label{alg:RepTwiConIdPolycyclic}
	\begin{algorithmic}[1]
		\Function{RepTwistConjToId}{$\varphi,\psi,g$}
		\State \Return \textsc{RepTwistConjToIdByNormal}\((\varphi,\psi,g,G')\)
		\EndFunction
	\end{algorithmic}
\end{algorithm}

\begin{proposition}
	\label{prop:recursivereidclasses}
	Let \(G\) be a polycyclic group and \(\varphi,\psi \in \End(G)\). Let \(\bvarphi,\bpsi\) be the induced endomorphisms on the abelianisation \(G/G'\). Then \(R(\varphi,\psi)\) is finite if and only if \(R(\bvarphi,\bpsi)\) is finite and \(R(\iota_g\varphi|_{G'},\psi|_{G'})\) is finite for every \(g \in G\).
\end{proposition}

\begin{proof}
	First assume that \(R(\varphi,\psi) < \infty\). By \Cref{cor:psiisbijection}, then \(R(\iota_g\varphi,\psi) < \infty\) for all \(g \in G\), and by applying \Cref{prop:recursiveReidNrsPolycyclic} we indeed find that \(R(\bvarphi,\bpsi) < \infty\) and \(R(\iota_g\varphi|_{G'},\psi|_{G'}) < \infty\) for every \(g \in G\).	Conversely, assume that \(R(\bvarphi,\bpsi) < \infty\) and \(R(\iota_g\varphi|_{G'},\psi|_{G'}) < \infty\) for every \(g \in G\). By \Cref{thm:reidclassesunion}, \(\R(\varphi,\psi)\) is then a finite union of finite sets and hence \(R(\varphi,\psi) < \infty\).
\end{proof}
 \Cref{alg:RepReidClassPolycyclic} provides an implementation of \textsc{RepsReidClasses} for polycyclic groups of derived length at least \(2\).

 \begin{theorem}
 	Let \(G\) be a polycyclic group of derived length at least \(2\) and let \(\varphi,\psi \in \End(G)\). Then \(\varphi\), \(\psi\) and \(G'\) satisfy the conditions necessary to apply \Cref{alg:RepReidClassNormal}.
 \end{theorem}
 \begin{proof} We prove this condition by condition.

 	\begin{enumerate}[(i)]
 		\item[(i)]- (iii) These can be proven in the same way as \Cref{thm:algworks}.
 		\item[(iv)] This follows from \Cref{prop:recursivereidclasses}.
 		\item[(v)] \Cref{alg:RepTwiConIdPolycyclic} provides this implementation. \qedhere
 	\end{enumerate}
 \end{proof}

\begin{algorithm}
	\caption{Determining representatives of \(\R(\varphi,\psi)\) if \(G\) is polycyclic}
	\label{alg:RepReidClassPolycyclic}
	\begin{algorithmic}[1]
		\Function{RepsReidClasses}{$\varphi,\psi$}
		\State \Return \textsc{RepsReidClassesByNormal}\((\varphi,\psi,G')\)
		\EndFunction
	\end{algorithmic}
\end{algorithm}

\section{Polycyclic-by-finite groups}

We can extend the algorithms \textsc{RepTwistConjToId} and \textsc{RepsReidClasses} to polycyclic-by-finite groups as in \Cref{alg:RepTwiConIdPolycyclicByFinite,alg:RepReidClassPolycyclicByFinite}. Unlike for abelian and polycyclic groups, these algorithms are not practical. In order to obtain practical algorithms, we would require the following:

\begin{itemize}
	\item practical algorithms that make polycyclic-by-finite groups suitable for computation with homomorphisms, in the sense described at the end of \Cref{sec:intro},
	\item a practical algorithm that finds a fully invariant, finite index, polycyclic subgroup \(N\) of a given polycyclic-by-finite group \(G\).
\end{itemize}

\begin{algorithm}
	\caption{Determining \(h\) such that \(g = \psi(h) \varphi(h)^{-1}\) if \(G\) is polycyclic-by-finite}
	\label{alg:RepTwiConIdPolycyclicByFinite}
	\begin{algorithmic}[1]
		\Function{RepTwistConjToId}{$\varphi,\psi,g$}
		\State \(N := \) fully invariant, finite index, polycyclic subgroup of \(G\)
		\State \Return \textsc{RepTwistConjToIdByNormal}\((\varphi,\psi,g,N)\)
		\EndFunction
	\end{algorithmic}
\end{algorithm}
\begin{theorem}
	\label{thm:algworks2}
	Let \(G\) be a polycyclic-by-finite group, let \(\varphi,\psi \in \End(G)\) such that  \(R(\varphi,\psi) < \infty\) and let \(N\) be a fully invariant, finite index, polycyclic subgroup of \(G\). Then \(\varphi\), \(\psi\) and \(N\) satisfy the conditions needed to apply \Cref{alg:RepTwiConIdNormal}.
\end{theorem}
\begin{proof}We prove this condition by condition. 
	\begin{enumerate}[(i)]
		\item This condition is satisfied because \(N\) is fully invariant.
		\item Since \(\Coin(\bvarphi,\bpsi)\) is a subgroup of the finite quotient \(G/N\), it is finite and can be computed by comparing \(\bvarphi(\bar{g})\) and \(\bpsi(\bar{g})\) for every \(\bar{g} \in G/N\).
		\item Since \(G/N\) is finite, we can implement \textsc{RepTwistConjToId} for the induced endomorphisms \(\bvarphi,\bpsi\), e.g. by computing \(\bpsi(\bar{h})\bvarphi(\bar{h})^{-1}\) for every \(\bar{h} \in G/N\) and comparing with the input \(\bar{g}\).
		\item \Cref{alg:RepTwiConIdPolycyclic} provides an implementation for endomorphisms of \(N\).\qedhere
	\end{enumerate}
\end{proof}

\begin{algorithm}
	\caption{Determining representatives of \(\R(\varphi,\psi)\) if \(G\) is polycyclic-by-finite}
	\label{alg:RepReidClassPolycyclicByFinite}
	\begin{algorithmic}[1]
		\Function{RepsReidClasses}{$\varphi,\psi$}
		\State \(N := \) fully invariant, finite index, polycyclic subgroup of \(G\)
		\State \Return \textsc{RepsReidClassesByNormal}\((\varphi,\psi,N)\)
		\EndFunction
	\end{algorithmic}
\end{algorithm}
 \begin{theorem}
	Let \(G\) be a polycyclic-by-finite group, let \(\varphi,\psi \in \End(G)\) such that  \(R(\varphi,\psi) < \infty\) and let \(N\) be a fully invariant, finite index, polycyclic subgroup of \(G\). Then \(\varphi\), \(\psi\) and \(N\) satisfy the conditions needed to apply \Cref{alg:RepReidClassNormal}.
\end{theorem}
\begin{proof} We prove this condition by condition.
	
	\begin{enumerate}[(i)]
		\item[(i)]- (iii) These can be proven in the same way as \Cref{thm:algworks2}.
		\item[(iv)] This follows from \Cref{cor:2props}(2).
		\item[(v)] \Cref{alg:RepTwiConIdPolycyclicByFinite} provides this implementation. \qedhere
	\end{enumerate}
\end{proof}

Significant progress towards making computation in polycyclic-by-finite feasible has been made by Sinanan and Holt \cite{sh17-1}, although their algorithms do not include computations with homomorphisms, and they require prior knowledge of a finite index, polycyclic, normal subgroup. For the latter requirement a theoretical algorithm exists, as proven in the following proposition, but it is not practical.

\begin{proposition}
	There is an algorithm which finds a fully invariant, finite index, polycyclic subgroup \(N\) of a polycyclic-by-finite group \(G\).
\end{proposition}
\begin{proof}
	There exists an algorithm to find a finite index, polycyclic, normal subgroup \(P\) of \(G\) (see \cite[Proposition 2.8]{bcr91-1}). Let \(m\) be the exponent of the finite quotient \(G/P\), i.e. the smallest positive integer such that \(\bar{g}^m = \bar{1}\) for any \(\bar{g} \in G/P\). Then \(N := \langle g^m \mid g \in G \rangle\) is a fully invariant, finite index, polycyclic subgroup of \(G\). There exists an algorithm to find such a subgroup of a  polycyclic-by-finite group (see \cite[Proposition 2.10]{bcr91-1}).
\end{proof}

 While \Cref{alg:RepTwiConIdPolycyclicByFinite,alg:RepReidClassPolycyclicByFinite} currently cannot be implemented in general, if one has a polycyclic-by-finite group \(G\), a representation of \(G\) suitable for computation (e.g. a particular matrix representation) and knowledge of a finite index, polycyclic, normal subgroup \(N\) invariant under the endomorphisms \(\varphi\) and \(\psi\), the algorithms can be implemented for that specific case. For example, in \cite{dkt17-2} Reidemeister numbers of the form \(R(\varphi,\id)\) with \(\varphi \in \Aut(G)\) were calculated for crystallographic groups \(G\). \Cref{alg:RepReidClassPolycyclicByFinite} reduces to \cite[Algorithm 3]{dkt17-2} if \(G\) is crystallographic, \(\varphi \in \Aut(G)\), \(\psi = \id\) and \(N = \Fitt(G)\), the Fitting subgroup of \(G\).
 
\section{Implementation in GAP}
\Cref{alg:RepTwistConj,alg:RepTwiConIdNormal,alg:RepReidClassNormal,alg:RepTwiConIdAbelian,alg:RepTwiConIdPolycyclic,alg:RepReidClassAbelian,alg:RepReidClassPolycyclic} have been implemented in the computer algebra system GAP \cite{gap18-3}, as part of a package called \texttt{TwistedConjugacy} \cite{gap20-1}. Below, we give a short demonstration of how to access our algorithms using this package. By way of example, let \(G\) be the group given by the following presentation:
\begin{equation*}
G := \left\langle g_1,g_2,g_3,g_4 \;\Bigg\vert \begin{array}{ll}
[g_1,g_2]	= g_2^2		& [g_1,g_4] 	= 1 \\\relax
[g_1,g_3] 	= g_3^2 	& [g_2,g_4] = 1\\\relax
[g_2,g_3]	= g_4^{-2} 		& [g_3,g_4] = 1\\
g_1^2		= g_4		& 
\end{array} \right\rangle.
\end{equation*}
This is a polycyclic group of derived length \(3\), and can be accessed in GAP through the command \texttt{ExamplesOfSomePcpGroups} provided by the \texttt{polycyclic} package \cite{gap18-5}. Let \(\varphi\) and \(\psi\) be the endomorphisms of \(G\) given by
\begin{align*}
\varphi(g_1)	&= g_1g_4^{-1},		&\psi(g_1)	&= g_1,\\
\varphi(g_2)	&= g_3,				&\psi(g_2)	&= g_2^2g_3g_4^2,\\
\varphi(g_3)	&= g_2g_3^3g_4^3,	&\psi(g_3)	&= g_2g_3g_4,\\
\varphi(g_4)	&= g_4^{-1},		&\psi(g_4)	&= g_4.
\end{align*}
One may load the \texttt{TwistedConjugacy} package and construct \(G\), \(\varphi\) and \(\psi\) as follows.
\begin{verbatim}
gap> LoadPackage("TwistedConjugacy");;
gap> G := ExamplesOfSomePcpGroups( 5 );;
gap> gens := GeneratorsOfGroup( G );;
gap> imgs1 := [ G.1*G.4^-1, G.3, G.2*G.3^3*G.4^3, G.4^-1  ];;
gap> phi := GroupHomomorphismByImages( G, G, gens, imgs1 );
[ g1, g2, g3, g4 ] -> [ g1*g4^-1, g3, g2*g3^3*g4^3, g4^-1 ]
gap> imgs2 := [ G.1, G.2^2*G.3*G.4^2, G.2*G.3*G.4, G.4  ];;
gap> psi := GroupHomomorphismByImages( G, G, gens, imgs2 );
[ g1, g2, g3, g4 ] -> [ g1, g2^2*g3*g4^2, g2*g3*g4, g4 ]
\end{verbatim}
The command \texttt{RepresentativeTwistedConjugation} provides an implementation of the \textsc{RepTwistConj} algorithm. We can use it to show that \(g_1\) and \(g_1^2\) are not \((\varphi,\psi)\)-twisted conjugate and that \(g_1\) and \(g_1^3\) are. Note that if two elements are not twisted conjugate, this implementation will return \texttt{fail} (rather than \AlgFail).
\begin{verbatim}
gap> RepresentativeTwistedConjugation( phi, psi, G.1, G.1^2 );
fail
gap> RepresentativeTwistedConjugation( phi, psi, G.1, G.1^3 );
g1*g4^-1 
\end{verbatim}
The command \texttt{ReidemeisterClasses} provides an implementation of the \textsc{RepsReidClasses} algorithm. We use it to show that \(R(\id,\psi) = \infty\) and to calculate representatives of the Reidemeister classes of \((\varphi,\psi)\).  Note that if the Reidemeister number is infinite, this implementation will return \texttt{fail} (rather than \AlgFail).
\begin{verbatim}
gap> ReidemeisterClasses( IdentityMapping( G ), psi );
fail
gap> ReidemeisterClasses( phi, psi );
[ id^G, g1*g2*g3^G, g1*g2^G, g1*g3^G, g1^G, g2*g3^G, g2^G, g3^G ] 
\end{verbatim}
Note that the ``\verb|^G|'' in the output above indicates that these elements are representatives of the orbits of a group action. For more information on the \texttt{TwistedConjugacy} package for GAP, we refer to \href{https://sTertooy.github.io/TwistedConjugacy/doc/chap0.html}{the package manual}.

\section*{Acknowledgement}
The authors would like to thank the referee for their careful reading and detailed comments.
\label{bibliography}
\addcontentsline{toc}{section}{References}
\printbibliography

\end{document}